\documentclass[10pt,letterpaper,conference]{ieeeconf}

\usepackage{latexsym, amssymb, amsmath}
\usepackage[dvips]{graphicx} 

\usepackage{float}      
\restylefloat{figure}

\usepackage{graphicx,dsfont}
\usepackage{mathrsfs}  
\usepackage{color} 
\usepackage[latin1]{inputenc}
\usepackage{enumerate}
\usepackage{flushend} 
\usepackage{stfloats}
\usepackage{cite} 

\usepackage{balance}

\allowdisplaybreaks[4] 

\newtheorem{prop}{Proposition} 
\newtheorem{proposition}[prop]{Proposition} 
\newtheorem{lem}{Lemma}
\newtheorem{lemma}[lem]{Lemma} 
\newtheorem{thm}{Theorem} 
\newtheorem{theorem}[thm]{Theorem} 
\newtheorem{cor}{Corollary} 
\newtheorem{corollary}[cor]{Corollary} 
\newtheorem{defn}{Definition} 
\newtheorem{definition}[defn]{Definition} 
\newtheorem{exmp}{Example} 
\newtheorem{example}[exmp]{Example} 
\newtheorem{exam*}{Example}

\def\custombibliography#1{
 \normalsize
\section*{\centering References}
 \list
 {[\arabic{enumi}]}{\settowidth\labelwidth{[#1]}\leftmargin\labelwidth
 \setlength{\itemsep}{.1em}
 \advance\leftmargin\labelsep
 \usecounter{enumi}}
 \def\newblock{\hskip .11em plus .33em minus -.07em}
 \sloppy
 \sfcode`\.=1000\relax}

\def\L2{{\cal L}_2}
\def\bull{\rule{0.08in}{0.08in}} 
\def\openbull{\framebox[0.08in][c]{$\;$}} 
\def\re{{\mathbb R}} 

\def\C{{\mathbb C}} 

\def\ss#1{{\scriptstyle #1}} 
\def\eqref#1{(\ref{#1})} 

\def\Tr{{\rm Tr}}
\def\vec{{\rm vec}}

\newcommand{\comment}[1]{} 

\def\begce{\begin{center}}
\def\endce{\end{center}}
\def\begar{\begin{array}}
\def\endar{\end{array}}
\def\begeq{\begin{equation}}
\def\endeq{\end{equation}}
\def\begdi{\begin{displaymath}}
\def\enddi{\end{displaymath}}
\def\begdis{\begin{eqnarray*}}
\def\enddis{\end{eqnarray*}}
\def\begeqa{\begin{eqnarray}}
\def\endeqa{\end{eqnarray}}
\def\begdes{\begin{description}}
\def\enddes{\end{description}}
\def\begit{\begin{itemize}}
\def\endit{\end{itemize}}
\def\begen{\begin{enumerate}}
\def\enden{\end{enumerate}}
\def\beglar{\left[\begin{array}}
\def\endrar{\end{array}\right]}
\def\begle{\begin{lemma}}
\def\endle{\end{lemma}}
\def\begde{\begin{definition}}
\def\endde{\end{definition}}
\def\begth{\begin{theorem}}
\def\endth{\end{theorem}}
\def\begco{\begin{corollary}}
\def\endco{\end{corollary}}
\def\begprop{\begin{proposition}}
\def\endprop{\end{proposition}}
\def\begex{\begin{example}}
\def\endex{\hfill\openbull \end{example} \vspace*{0.1in}}
\def\begexer{\begin{exercise}}
\def\endexer{\end{exercise}}

\def\begres{\noindent{\bf Remarks}:\begin{enumerate}}
\def\endres{\end{enumerate} \par}
\def\begpr{\noindent{\em Proof:}$\;\;$}
\def\endpr{\hfill\bull \vspace*{0.1in}}

\def\begtab{\begin{tabular}}
\def\endtab{\end{tabular}}
\def\rref#1{(\ref{#1})}

\newcommand\cdcout[1]{} 

\newcommand{\rv}[1]{\boldsymbol{#1}} 
\newcommand{\RomanNumber}[1]{\uppercase\expandafter{\romannumeral #1}}
\newcommand{\romannumber}[1]{\lowercase\expandafter{\romannumeral #1}}

\DeclareMathAlphabet{\mathpzc}{OT1}{pzc}{m}{it}

\def\1{\rv 1} 

\def\allseriesA{\mbox{$\re\langle\langle A\rangle\rangle$}}
\def\allseriesA'{\mbox{$\re\langle\langle A'\rangle\rangle$}}

\def\L1spaceprodu{{ L}_1(\Omega\times [0,T],{\mathcal P},P\otimes \lambda)}

\def\Hspace0{{\mathcal H}^2_0}

\newcount\colveccount
\newcommand*\colvec[1]{
        \global\colveccount#1
        \begin{pmatrix}
        \colvecnext
}
\def\colvecnext#1{
        #1
        \global\advance\colveccount-1
        \ifnum\colveccount>0
                \\
                \expandafter\colvecnext
        \else
                \end{pmatrix}
        \fi
}

\renewcommand{\baselinestretch}{0.95}

\usepackage{cite} 
\usepackage{color} 
\usepackage{graphicx} 
\usepackage{latexsym,amssymb,amsmath} 
\usepackage{mathrsfs} 
\usepackage{multicol}
\usepackage[german,english]{babel}
\usepackage[T1]{fontenc}
\usepackage[latin1]{inputenc}

\allowdisplaybreaks[4]

\IEEEoverridecommandlockouts                              
\title{On the Preservation of Commutation and Anticommutation Relations of $n$-Level Quantum Systems$^\ast$\thanks{$^\ast$This work was supported by the
Australian Research Council (ARC) projects FL110100020 and DP110102322, and Air Force Office of Scientific Research (AFOSR). This material is based on research
sponsored by the Air Force Research Laboratory, under agreement number FA2386-09-1-4089. The U.S. Government is authorized to reproduce and distribute reprints
for Governmental purposes notwithstanding any copyright notation thereon. The views and conclusions contained herein are those of the authors and should not be
interpreted as necessarily representing the official policies or endorsements, either expressed or implied, of the Air Force Research Laboratory or the U.S.
Government. }}

\author{Luis~A.~Duffaut~Espinosa$^\dagger$\thanks{$^\dagger$School of Engineering and Information Technology, University of New South Wales at ADFA, Canberra,
ACT 2600, Australia. {\tt\small \{l.duffaut,  i.petersen, v.ougrinovski\}@adfa.edu.au.}}, Z.~Miao$^\ddag$\thanks{$^\ddag$Research School of Engineering, The
Australian National University, Canberra ACT 0200, Australia. {\tt\small zibo.miao@anu.edu.au}.}, I.~R.~Petersen$^\dagger$, V.~Ugrinovskii$\,^{\dagger
\S}$\thanks{$\,^{\S}$Part of this work was carried out during this author visit to the Australian National University}, and M.~R.~James$^\$$\thanks{$^\$$ARC
Centre for Quantum Computation and Communication Technology, Research School of Engineering, Australian National University, Canberra, ACT 0200, Australia.
{\tt\small matthew.james@anu.edu.au}.}}

\date{\today}

\begin{document}

\maketitle

\begin{abstract}
The goal of this paper is to provide conditions under which a quantum stochastic differential equation (QSDE) preserves the commutation and anticommutation
relations of the $\pmb{SU(n)}$ algebra, and thus describes the evolution of an open $\pmb{n}$-level quantum system. One of the challenges in the approach lies
in the handling of the so-called anomaly coefficients of $\pmb{SU(n)}$. Then, it is shown that the physical realizability conditions recently developed by the
authors for open $\pmb{n}$-level quantum systems also imply preservation of commutation and anticommutation relations. 
\end{abstract}

\section{Introduction}  \label{sec:section1}

The property of \emph{physical realizability} of open quantum systems has attracted considerable interest in recent years with its main motivation being its
role in {\em quantum coherent control} \cite{Bouten-Handel-James_2007,Sarovar-Ahn-Jacobs-Milburn_2004,Helon-James_2006,Lloyd_2000}. In simple words,
this property systematically characterizes the quantum nature of a control system from a state space point of view, which is relevant in many engineering areas,
and is closely related to algebraic properties of the QSDE coefficients. The case of open quantum linear systems is well understood
\cite{James-Nurdin-Petersen_2008}. More recently, the physical realizability property was introduced for $n$-level quantum systems \cite{Duffaut-et-al_2012c}.
This was the product of a nontrivial extension of the two-level case presented in \cite{Duffaut-et-al_2012a}. From the state space point of view, physical
realizability has to be complemented by the preservation of commutation relations since the latter is an exhibitor of the quantum behavior of the system. In
this regard, the authors provided conditions for the preservation of commutation relations in the case of QSDE's evolving in $SU(2)$ \cite{Duffaut-et-al_2012b}.
However, in order to extend the formalism to the general $n$-level case it is necessary to deal with the so-called \emph{anomaly coefficients of} $SU(n)$, which
forms the completely symmetric tensor $d_{ijk}$ that appears out of the anticommutation relation of the generators of $SU(n)$. These generators form a complete
orthonormal basis spanning the set of $n$-dimensional complex matrices and are known as the {\em generalized Gell-Man matrices}
\cite{Mahler-Weberrus_98,Macfarlane-Sudbery-Weisz_68,Kaplan-Resnikoff_67}. The commutation relations of these generators also generate the completely
antisymmetric tensor $f_{ijk}$. In the case of QSDEs describing two-level quantum systems evolving in $SU(2)$, it suffices to consider preserving the
commutation relations of $SU(2)$ since the anomaly tensor $d$ is zero for all indices. On the other hand, the situation is nontrivial for $n\ge 3$ since one now
requires extra machinery on the tensor $d$ that in some cases depend on the value of $n$. In other words, the QSDE for an $n$-level quantum system must preserve
both commutation and anticommutation relations for all times. 

The approach followed in this manuscript is that of open quantum systems defined in terms of a triple $(S,L,{H})$, where $S$ is an operator, known as the {\em
scattering matrix}, describing the interaction of the environment fields among themselves, $L$ is a vector of coupling operators expressing the interaction of
the environment fields with the system variables and $H$ is a Hamiltonian operator \cite{Gough-James_2009,James-Gough_2010}. This description implicitly
takes in consideration the quantum nature of the evolution equations. For example, the $(S,L,H)$ description will preserve, in time, the canonical commutation
relations of the system under consideration. The physical realizability conditions in \cite{Duffaut-et-al_2012c} do {\em not}, in general, assure the
preservation of commutation and anticommutation relations. It is thus the main goal of this paper to provide conditions for the preservation of commutation and
anticommutation relations for a state space model evolving in $SU(n)$ independently of the physical realizability conditions. Then, as a second result it is
shown that physical realizability \emph{does} imply the preservation of commutation and anticommutation relations.

The paper is organized as follows. Section \ref{sec:section2} summarizes the necessary tools obtained from the algebra of $SU(n)$. In Section \ref{sec:section3}
the basic preliminaries of open $n$-level quantum systems as well as the property of physical realizability for such systems are introduced. This is followed by
Section \ref{sec:section4}, in which conditions for the preservation of commutation and anticommutation relations are provided. In addition, it is also shown
that those conditions are implied by the physical realizability conditions. Finally, Section \ref{sec:conclusions} gives the conclusions.

\section{The algebra of $SU(n)$} \label{sec:section2}

In what follows the necessary tools regarding the algebra of the special unitary group $SU(n)$ will be provided; see
\cite{Mahler-Weberrus_98,Macfarlane-Sudbery-Weisz_68,Kaplan-Resnikoff_67} for more details. Basically, this group is formed by all $n\times n $ complex matrices
that are Hermitian and have zero trace. Consider the set of elementary vectors spanning $\C^n$, namely, $\{e_1, \cdots, e_n\}$. Define ${P}_{kl}\in \C^{n\times
n}$ as ${P}_{k,l} = e_k e_{l}^T$, where $k,l=1,\hdots,n$. A standard way of constructing a complete orthonormal basis for $SU(n)$ is 
\begin{align*}
u_{jk} & = P_{j,k}+P_{k,j},\\  
v_{jk} & = \pmb{i}\left(P_{j,k}-P_{k,j} \right),\\
w_l    & = -\sqrt{\frac{2}{l(l+1)}} \left( \sum_{s=1}^k P_{s,s} - k P_{l+1,l+1} \right) 
\end{align*}
for $1\le j < k \le n$, $1 \le l \le n-1$. Note that the identity matrix $I$ must be included in the basis in order to form a complete set. The identity,
$(n^2-n)/2$ symmetric matrices $u_{jk}$, the $(n^2-n)/2$ antisymmetric matrices $v_{jk}$ and the $n-1$ mutually commutative matrices $w_l$ together form the
\emph{generators of $SU(n)$}. These generators are known as the \emph{generalized Gell-Mann matrices}. Without any particular order, the generators are
relabeled $\{I,\lambda_1, \hdots, \lambda_{s} \}$, where $s = n^2-1$. Here the orthonormality condition these matrices satisfy is $\Tr(\lambda_i \lambda_j) =2
\delta_{ij}$, where $\delta_{ij}$ denotes the Kronecker delta. Their commutation and anticommutation relations are
\begin{align*}
[\lambda_i,\lambda_j] & = 2 \pmb{i} \sum_{k=1}^{s} f_{ijk} \lambda_{k}, \\
\{\lambda_i,\lambda_j \} & = \frac{4}{n} \delta_{ij} + 2 \sum_{k=1}^{s} d_{ijk} \lambda_{k}. 
\end{align*}
Thus, the product $\lambda_i \lambda_j$ can be easily computed as
\begin{align} \label{eq:product_Gell-mann_matrices}
\nonumber \lambda_i \lambda_j & = \frac{1}{2} \left( [\lambda_i,\lambda_j] + \{ \lambda_i,\lambda_j \} \right) \\ 
                    & = \frac{2}{n} \delta_{ij} + \sum_{k=1}^{s} \left( \pmb{i} f_{ijk} +d_{ijk}  \right)\lambda_k.
\end{align}
where the real completely antisymmetric tensor $f_{ijk}$ and the real completely symmetric tensor $d_{ijk}$ are called the \emph{structure constants} of
$SU(n)$. The tensors $f_{ijk}$ and $d_{ijk}$ satisfy 
\begin{subequations}
\begin{alignat}{2} 
\label{eq:property_ff_tensor} & f_{ilm}f_{mjk}+f_{jlm}f_{imk}+f_{klm}f_{ijm}=0, \\
\label{eq:property_fd_tensor} & f_{ilm}d_{mjk}+f_{jlm}d_{imk}+f_{klm}d_{ijm}=0, \\
\nonumber & \sum_{k=1}^s f_{ilk}f_{mjk} =  \frac{2}{n}\left(\delta_{im} \delta_{lj} - \delta_{ij} \delta_{lm} \right) \\ 
\label{eq:property_ffk_tensor}                                        & \hspace{0.9in} + \sum_{k=1}^s \left( d_{imk} d_{ljk} - d_{ijk} d_{lmk} \right), \\
\label{eq:property_ff_delta}  &\sum_{m,k=1}^s f_{im k} f_{jm k}  = n  \delta_{ij}. 
\end{alignat}
\end{subequations}
Define $F_i,D_i \in \re^{s\times s}$, $i\in \{1,\hdots, s\}$, such that their $(j,k)$ component is $({F}_i)_{jk} = f_{ijk}$ and $(D_i)_{jk} = d_{ijk}$,
respectively. In particular, the set $\{-\pmb{i} {F}_1, \hdots, -\pmb{i} {F}_{s}  \}$ is the adjoint representation of $SU(n)$. In
\cite{Macfarlane-Sudbery-Weisz_68,Kaplan-Resnikoff_67}, identities \rref{eq:property_ff_tensor}-\rref{eq:property_ffk_tensor} were employed to obtain the
following useful relationships
\begin{subequations}
    \begin{alignat}{2}
\label{eq:D_F_3}      [F_i,F_j] = & - \sum_{k}^{s} f_{ijk} F_k \\
\label{eq:D_F_4}      [F_i,D_j] = & - \sum_{k}^{s} f_{ijk} D_k \\
\label{eq:D_F_5}      F_i D_j + F_jD_i = &  \sum_{k}^{s} d_{ijk} F_k \\
\label{eq:D_F_6}      D_i F_j + D_jF_i = &  \sum_{k}^{s} d_{ijk} F_k\\
\nonumber             \left(D_iD_j - F_jF_i\right)_{ml} = & \sum_{k}^{s} d_{ijk} (D_k)_{ml} \\ 
\label{eq:D_D_7}                                       & \, + \frac{2}{n} \left(\delta_{ij}\delta_{ml}-\delta_{im}\delta_{jl}\right).
    \end{alignat}
\end{subequations}

\begde
Let $\beta\in \C^s$. The linear mappings $\Theta^-,\Theta^+: \C^s \rightarrow
\C^{s\times s}$ are defined as 
\begin{align*}
\Theta^-(\beta) &= \left(F_1^T \beta,\cdots , F_{s}^T \beta \right) = \left(\begin{array}{c}
         \beta^T F_1^T  \\ \vdots \\ \beta^T F_{s}^T  \end{array} \right),\\
\Theta^+(\beta) & = \left(D_1^T \beta,\cdots , D_{s}^T \beta \right) = \left(\begin{array}{c}
\beta^T D_1^T  \\ \vdots \\ \beta^T D_{s}^T  \end{array} \right).
\end{align*}  
\endde
\vspace*{0.05in}

\noindent Observe that the nature of the $f$ and $d$-tensors make $\Theta^-(\beta)$ and $\Theta^+(\beta)$ be antisymmetric and symmetric, respectively. When
$\beta$ is an $s$-dimensional row vector then it will be understood hereafter that $\Theta^-(\beta) = \Theta^-(\beta^T)$ and
$\Theta^+(\beta)=\Theta^+(\beta^T)$. Consider now the \emph{stacking operator} $\vec: \C^{m\times n} \rightarrow \C^{m n}$ whose action on a matrix creates a
column vector by stacking its columns below one another. With the help of $\vec$, the matrices $\Theta^-(\beta)$ and $\Theta^+(\beta)$ can be reorganized so
that 
\begdi
\vec(\Theta^-(\beta))=\left(\begin{array}{c} \Theta^-_1(\beta) \\ \vdots \\ \Theta^-_s(\beta) \end{array} \right)= F \beta, 
\enddi
and
\begdi
\vec(\Theta^+(\beta))=\left(\begin{array}{c} \Theta^+_1(\beta) \\ \vdots \\ \Theta^+_s(\beta) \end{array} \right)= D \beta, 
\enddi
where $\beta\in \C^s$, $\Theta^-_i(\beta)=F^T_i\beta$, $F=\left(F_1 ,\cdots ,F_s  \right)^T$, $\Theta^+_i(\beta)=D_i\beta$ and $D=\left(D_1 ,\cdots,D_s
\right)^T$. From \rref{eq:property_ff_delta}, $F$ satisfies
\begeq \label{eq:E^T_by_E}
F^T  F= n I.
\endeq

The properties of $\Theta^-$ and $\Theta^+$ are summarized in the next lemma. 
\begle \label{le:prop_thetas} Let $\beta,\gamma\in \C^s$. The mappings $\Theta^-$ and $\Theta^+$ satisfy
 \begin{subequations}
    \begin{alignat}{2}
\label{eq:Theta_1} &  \hspace{-0.1in}\Theta^-(\beta) \gamma = -\Theta^-(\gamma) \beta, \\
\label{eq:Theta_2} & \hspace{-0.1in}\Theta^+(\beta) \gamma = \Theta^+(\gamma) \beta,\\
\label{eq:Theta_beta_beta} & \hspace{-0.1in}\Theta^-(\beta)\beta = 0, \\
\label{eq:Theta__minus_theta_minus} & \hspace{-0.1in}\Theta^-\left(\Theta^-(\beta) \gamma \right) = [\Theta^-(\beta) ,\Theta^-(\gamma) ],  \\
\label{eq:Theta__minus_theta_plus} & \hspace{-0.1in}\Theta^-\left(\Theta^+(\beta) \gamma \right) = \Theta^-(\beta) \Theta^+(\gamma) + \Theta^-(\gamma)
\Theta^+(\beta), \\
\label{eq:Theta__plus_theta_minus1}    & \hspace{-0.1in}\Theta^+\left(\Theta^-(\beta) \gamma \right) = [\Theta^+(\beta),\Theta^-(\gamma)] = [\Theta^-(\beta),
\Theta^+(\gamma) ],\\
\nonumber                      & \hspace{-0.1in} \Theta^+\left( \Theta^+(\beta) \gamma  \right) = \Theta^+(\beta) \Theta^+(\gamma) -
\Theta^-(\gamma)\Theta^-(\beta) \\
\label{eq:thetaplus_thetaplus} & \hspace{-0.1in}\hspace{1.0in} - \frac{2}{n}\left( \beta^T \gamma I - \beta \gamma^T\right).
    \end{alignat}
\end{subequations}
\endle\vspace*{0.05in}

\begpr The proofs of identities \rref{eq:Theta_1}-\rref{eq:Theta__plus_theta_minus1} where provided by the authors in \cite{Duffaut-et-al_2012c}. Only
\rref{eq:thetaplus_thetaplus} requires a proof. The left-hand-side of \rref{eq:thetaplus_thetaplus} is decomposed as 
\begin{align*}
\lefteqn{\Theta^+\left(\Theta^+(\beta) \gamma \right)} \\
 & = \left(  D_1 \left(\begin{array}{c} \beta^T D_1 \gamma \\ \vdots \\  \beta^T D_s \gamma \end{array}\right), \cdots, D_s \left(\begin{array}{c} \beta^T D_1
\gamma \\ \vdots \\  \beta^T D_s \gamma \end{array}\right) \right) \\
& = \left( \begin{array}{ccc} \displaystyle \sum_{k=1}^{s} d_{11k} \beta^T D_k \gamma & \cdots & \displaystyle \sum_{k=1}^{s} d_{s1k} \beta^T D_k \gamma  \\ 
\vdots &  \ddots &  \vdots  \\ \displaystyle \sum_{k=1}^{s} d_{1sk} \beta^T D_k \gamma & \cdots & \displaystyle  \sum_{k=1}^{s} d_{ssk} \beta^T D_k \gamma 
\end{array} \right).
\end{align*}
By \rref{eq:D_D_7}, the $(i,j)$ component of this matrix is 
\begin{align*}
\lefteqn{\hspace*{-0.3in} \left( \Theta^+\left(\Theta^+(\beta) \gamma \right) \right)_{ij} } \\
 = & \, \sum_{k=1}^{s} d_{ijk} \beta^T D_k \gamma \\
 = & \, \beta^T \left(\sum_{k=1}^{s} d_{ijk} D_k \right) \gamma \\
 = & \, \beta^T \left(D_iD_j - F_jF_i\right) \gamma \\
   &  - \frac{2}{n} \left(\delta_{ij}\sum_{m,l=1}^s\beta_m \delta_{ml} \gamma_l - \sum_{m,l=1}^s\beta_m \delta_{im}\delta_{jl} \gamma_l \right)\\
 = & \, \beta^T D_iD_j \gamma - \beta^T F_jF_i \gamma  - \frac{2}{n} \left( \beta^T \gamma \delta_{ij} - \beta_i \gamma_j \right)\\
 = & \, \beta^T D_iD_j \gamma - \gamma^T F_i^TF_j^T \beta - \frac{2}{n} \left( \beta^T \gamma \delta_{ij} - \beta_i \gamma_j \right) \\
 = & \, \left( \Theta^+(\beta) \Theta^+(\gamma) - \Theta^-(\gamma) \Theta^-(\beta) \right)_{ij} \\
   &  - \frac{2}{n} \left( \beta^T \gamma I - \beta \gamma^T \right)_{ij}.\end{align*}
\endpr

Some additional identities regarding matrices $F$ and $D$ with respect to the Kronecker product are given next. They will become useful when proving the main
results of the paper in Section \ref{sec:section4}. Define the tensor permutation matrix $\mathds{1}_\otimes \in \re^{s^2\times s^2}$. That is, a
symmetric block matrix $\mathds{1}_\otimes = \{\mathds{1}_{ji}\}_{i,j=1}^s$ such that it satisfies $\mathds{1}_\otimes (A\otimes B) \mathds{1}_\otimes =
(B\otimes A)$ for any $A,B\in \C^{s\times s}$.
\begle \label{le:F_kronecker_properties}
Let $F=\left(F_1 , F_2, \cdots, F_s  \right)^T$ and $F=\left(D_1 , D_2, \cdots, D_s  \right)^T$, and $A,B\in \re^{s\times s}$. Then

\renewcommand*\theenumi{{\roman{enumi}}}
\renewcommand*\labelenumi{$\theenumi.$}
\begin{enumerate}
\item \label{eq:kron_permutation1} $F = -\mathds{1}_\otimes F,$ \vspace*{0.05in}
\item \label{eq:kron_permutation2} $D = \mathds{1}_\otimes D,$ \vspace*{0.05in}
\item \label{eq:kron_permutation3} $F^T(A\otimes B)F = F^T(B\otimes A)F,$ \vspace*{0.05in}
\item \label{eq:kron_permutation4} $D^T(A\otimes B)D = D^T(B\otimes A)D.$ \vspace*{0.05in}
\end{enumerate}
\endle 
\begpr Identity $\rref{eq:kron_permutation1}$ is proved by multiplying directly $\mathds{1}_\otimes $ and $F$. Therefore,
\begin{align*}
\mathds{1}_\otimes F & = - \left(\begin{array}{ccc}  \mathds{1}_{11} & \cdots & \mathds{1}_{s1}  \\ 
                                                \vdots  & \ddots & \vdots \\
                                              \mathds{1}_{1s} & \cdots & \mathds{1}_{ss}  \end{array} \right) \colvec{3}{F_1}{\vdots}{F_s} \\
                     & = - \left(\begin{array}{c} \colvec{3}{(F_1)_1}{ \vdots }{(F_s)_1}  \\  \vdots \\  \colvec{3}{(F_1)_s}{\vdots}{(F_s)_s} \end{array}
\right).
\end{align*}
Observe that $(F_i)_j = (f_{ij1}, \cdots, f_{ijs}) = - (F_i)_j$ since $f_{ijk}$ is an antisymmetric tensor. Using this fact in the previous expression
gives 
\begin{align*}
\mathds{1}_\otimes F & = -\left(\begin{array}{ccc}  F^T_1 \\
                                                \vdots \\
                                                   F^T_s  \\   \end{array} \right) = -F
\end{align*}
An analogous procedure is used to show identity $\rref{eq:kron_permutation2}$, but in this case the permutation of indices in $d_{ijk}$ does not produce a
negative sign due to the fact that $d_{ijk}$ is a completely symmetric tensor. For $\rref{eq:kron_permutation3}$, one will employ $\rref{eq:kron_permutation1}$
as follows
\begin{align*}
F^T(A\otimes B)F & = (-\mathds{1}_\otimes F)^T(A\otimes B)(-\mathds{1}_\otimes F)\\
		 & = F^T (\mathds{1}_\otimes (A\otimes B)\mathds{1}_\otimes) F \\
		 & = F^T(B\otimes A)F.
\end{align*}
Finally, identity $\rref{eq:kron_permutation4}$ is proved similarly but using $\rref{eq:kron_permutation2}$ instead. 
\endpr

\section{Physical Realizability of Open $n$-Level Quantum Systems} \label{sec:section3}

Consider the separable Hilbert space $\mathfrak{h}$, $\mathfrak{T}(\mathfrak{H})$ the set of operators in $\mathfrak{H}$, and
$\mathfrak{T}(\mathfrak{H})^{n\times m}$ the set of $n\times m$ dimensional arrays of operators in $\mathfrak{T}(\mathfrak{H})$. The commutator of $x$ and $y$
in $\mathfrak{T}({\mathfrak{H}})$ is $[x, y] = xy - xy$. Let $x\in \mathfrak{T}({\mathfrak{H}})^{n_1}$ and $y\in \mathfrak{T}({\mathfrak{H}})^{n_2}$, then 
$[{x},{y}^T] \triangleq {x} {y}^T - ({y} {x}^T)^T\in \mathfrak{T}({\mathfrak{H}})^{n_1\times n_2}$. 
The adjoint of $x$ is denoted by $x^\dagger = (x^\#)^T$ with 
\begdi 
{x}^\# \triangleq \colvec{3}{x_1^\ast}{ \vdots }{x_n^\ast} 
\enddi
and $^\ast$ denotes the operator adjoint. For the case of complex vectors and matrices, $^\ast$ denotes the complex conjugate while $^\dagger$ denotes the
conjugate transpose. 

The \emph{Heisenberg evolution} of $x\in \mathfrak{T}({\mathfrak{H}})^{s}$ interacting with $n_w$ quadrature Boson quantum fields
$\bar{W}_1$ and $\bar{W}_2$ in the parametrization $(S,L,H)$ is given by  
\begin{align} \label{eq:general_evolution_vector}
\nonumber dx = & \, \mathcal{L}(x)\,dt + \frac{1}{2}\left([x,L^T]-[x,L^\dagger] \right)d\bar{W}_1  \\ 
& \,  -\frac{\pmb{i}}{2}\left([x,L^T]+[x,L^\dagger] \right)d\bar{W}_2
\end{align}
where ${\cal L}(X)$ is the \emph{Lindblad operator} defined as
\begin{align}  \label{eq:Lindblad_operator_vector}
\nonumber \lefteqn{\hspace*{-0.2in} \mathcal{L}(x) } \\ 
            = & \, -\pmb{i} [x,H] + \frac{1}{2} \left( \left(L^\dagger \,  [x, L^T]^T   \right)^T   +[ L^\#,x^T]^T\,L  \right).
\end{align}
The output field in its quadrature form is
\begdi
\colvec{2}{d\bar{Y}_1}{d\bar{Y}_2} =\colvec{2}{L+L^\#}{\pmb{i}(L^\#-L) } \,dt + \colvec{2}{d\bar{W}_1}{d\bar{W}_2}.
\enddi
Here the operator $S$ is assumed to be the identity operator ($S=\hat{I}$), and the It\^o table (\cite{Hudson_parthasarathy_84}) for the fields $\bar{W}_1$ and
$\bar{W}_2$ is
\begeq \label{eq:Ito_table_quadrature}
\colvec{2}{d\bar{W}_1}{d\bar{W}_2} \left(\begin{array}{cc}
d\bar{W}_1 & d\bar{W}_2
\end{array}\right)= \left(\begin{array}{cc}
I_{n_w} & \pmb{i}I_{n_w} \\ -\pmb{i}I_{n_w} & I_{n_w} 
\end{array}\right) dt.
\endeq

The interest in this paper is in systems evolving with respect to the special unitary group $SU(n)$. Consider the Hilbert space for these systems to be
$\mathfrak{H}=\C^n$. It is standard to associate a vector $\beta=(\beta_1,\cdots, \beta_s)\in \C^s$ with the vector of operators $\hat{\beta} \in
\mathfrak{T}(\mathfrak{H})^s$ by simply considering the identity operator in each component of the vector, i.e., 
$\hat{\beta}=(\beta_1 \hat{I}, \cdots,\beta_s \hat{I})$. 
In a similar manner, any complex matrix is associated to a matrix of operators by considering the identity operator $\hat{I}$ in each component. The identity
operator $\hat{I}$ is usually suppressed. Therefore, abusing of the notation slightly, the mappings $\Theta^-$ and $\Theta^+$ are allowed to act on a vector of
operators. In this manner, it is valid to multiply complex matrices and operator matrices. This convention allows to express the commutation and
anticommutation relations of the system variables $x$ as 
\begin{subequations}
    \begin{align}
\label{eq:lambda_CCR}  [x,x^T]    & = 2 \pmb{i} \Theta^-(x),\\
\label{eq:lambda_antiCCR} \{x,x^T \} & = \frac{4}{n} I + 2 \Theta^+(x).
    \end{align}
\end{subequations}
The vector of system variables for \rref{eq:general_evolution_vector} evolving in $SU(n)$ is then
\begdi 
x=\colvec{3}{x_1}{ \vdots }{x_s} \triangleq \colvec{3}{\hat{\lambda}_1}{ \vdots }{\hat{\lambda}_s},
\enddi 
where $\hat{\lambda}_1, \ldots, \hat{\lambda}_s$ are self-adjoint and spanned by the generalized Gell-Mann matrices. They are usually called \emph{spin
operators}. The initial value of the system variables can be set to $x(0)=(\lambda_1,\ldots,\lambda_s)$ with $\lambda_1,\hdots, \lambda_s$ being the generators
of $SU(n)$ introduced in Section \ref{sec:section2}. Due to the product relation \rref{eq:product_Gell-mann_matrices} any polynomial of spin operators can be
written as a linear combination of generalized Gell-Mann matrices. Therefore, assuming linearity captures a large class of Hamiltonian and coupling operators. 
This allows to assume, without losing generality, that the Hamiltonian is ${\mathcal{H}}=\alpha x$ with $\alpha\in \re^{s}$, and the multiplicative coupling
operator is of the form $L=\Lambda x$ with $\Lambda \in \C^{n_w \times s}$.
    
In general, the evolution of $x$ in quadrature form falls into a class of bilinear QSDEs expressed as
\begin{align} \label{eq:bilinear_system}
\nonumber \lefteqn{ \hspace*{-0.1in} dx = A_0\,dt+A x \, dt} \\
& \hspace*{-0.1in}\,+ \left( B_{11}x , \cdots , B_{1n_w}x , B_{21}x  , \cdots , B_{2n_w} x \right)\colvec{2}{d\bar{W}_1}{d\bar{W_2}},
\end{align}
\begeq \label{eq:bilinear_system_output}
\colvec{2}{d\bar{Y}_1}{d\bar{Y}_2} = \colvec{2}{C_1}{C_2} x \,dt + \colvec{2}{d\bar{W}_1}{d\bar{W}_2},
\endeq
where $A_0\in \re^s$, $A,{B}_{1k}\triangleq \bar{B}_{1k}+\bar{B}_{2k},{B}_{2k}\triangleq \pmb{i}(\bar{B}_{2k}-\bar{B}_{1k}) \in \re^{s\times s}$ and $C_1,C_2\in
\re^{n_w \times s}$,
$k=1,\hdots, n_w$. It is assumed that all matrices in \rref{eq:bilinear_system} and \rref{eq:bilinear_system_output} are real due to the fact that the
class of quantum systems considered in this paper are in quadrature form. The next results in physical realizability were given in \cite{Duffaut-et-al_2012c}
and are employed later in Section \ref{sec:section4}.

\begde \cite{Duffaut-et-al_2012c} \label{def:physical realizability} A system described by equations \rref{eq:bilinear_system} and
\rref{eq:bilinear_system_output} is said to be {physically realizable} if there exist $\mathcal{H}$ and $L$ such that \rref{eq:bilinear_system} can be written
as in \rref{eq:general_evolution_vector}.
\endde

The explicit form of matrices $A_0,A, B_{1k}, B_{2k},C_1$ and $C_2$ in terms of the Hamiltonian and coupling operator is given next.
\begth \cite{Duffaut-et-al_2012c} \label{th:physical_realizability_def} Let $\mathcal{H}=\alpha x$, with $\alpha^T \in \re^s$, and $L=\Lambda x$, with $\Lambda
\in \C^{n_w \times s}$. Then
\begin{subequations}\label{eqn:Spin_system_matrices}
    \begin{align}
      A_0 &= \frac{4 \pmb{i}}{n} \sum_{k=1}^{n_w}\Theta^-(\Lambda^\#_k)\Lambda_k^T, \label{subeqn:F0}\\
      A   &= - 2 \Theta^-(\alpha) + \sum_{k=1}^{n_w} \left(R_k -\pmb{i} Q_k \right), \label{subeqn:F}\\
      {B}_{1k} &= \Theta^-\left({\pmb i}(\Lambda^\#_k - \Lambda_k)\right), \label{subeqn:G1}\\ 
      {B}_{2k} &= - \Theta^-(\Lambda_k + \Lambda^\#_k), \label{subeqn:G2} \\
      C_1   &=  \Lambda+\Lambda^\#, \label{subeqn:H12} \\
      C_2   &=  \pmb{i} \left( \Lambda^\#-\Lambda \right), \label{subeqn:H2}
   \end{align}
  \end{subequations}
where 
\begin{align*} 
R_k & \triangleq \Theta^-(\Lambda_k)\Theta^-(\Lambda^\#_k) + \Theta^-(\Lambda^\#_k)\Theta^-(\Lambda_k), \\
Q_k & \triangleq \Theta^-(\Lambda_k)\Theta^+(\Lambda^\#_k) - \Theta^-(\Lambda^\#_k)\Theta^+(\Lambda_k). 
\end{align*}
\endth \vspace{0.1in}

The next theorem present conditions for physical realizability in terms of $(A_0,A,B_i,C)$ in \rref{eq:bilinear_system} and
\rref{eq:bilinear_system_output}.

\begth \cite{Duffaut-et-al_2012c} \label{th:physical_realizability}
System \rref{eq:bilinear_system} with output equation \rref{eq:bilinear_system_output} is physically realizable if and only if 
\renewcommand*\theenumi{{\roman{enumi}}}
\renewcommand*\labelenumi{$\theenumi.$}
\begin{enumerate}
\item \label{itm:theorem_physical_realizability1} $\displaystyle A_0=\frac{1}{n}\sum_{k=1}^{n_w}(\pmb{i}{B}_{1k} + {B}_{2k}) \left((C_1)_k+\pmb{i}
(C_2)_k\right)^T$,
\vspace*{0.05in}
\item \label{itm:theorem_physical_realizability2} $\displaystyle {B}_{1k}= \Theta^-((C_2)_k)$,
\vspace*{0.1in}
\item \label{itm:theorem_physical_realizability21} $\displaystyle {B}_{2k}= \Theta^-((C_1)_k)$,
\vspace*{0.05in}
\item \label{itm:theorem_physical_realizability3} $\displaystyle A+A^T+\sum_{i,k=1}^{2,n_w} {B}_{ik} {{B}_{ik}}^T =\frac{n}{2}\Theta^+(A_0)$,
\vspace*{0.05in}
\end{enumerate}
where $(C_i)_k$ indicates the $k$-th row of $C_i$. In which case, the coupling matrix can be identified to be 
\begdi 
\Lambda=\frac{1}{2}(C_1+\pmb{i}C_2),
\enddi
and $\alpha$, defining the system Hamiltonian, is
\begin{align} \label{eq:Hamiltonian_physical_realizability}
\nonumber  \lefteqn{\hspace*{-0.8in}\alpha = \frac{1}{4n}\vec\left(A^T-A +\frac{1}{2} \sum_{k=1}^{n_w}\left(  [B_{2k} , \Theta^+((C_2)_k)] \right. \right. } \\
 & \, \left. \left. \rule{0in}{0.2in} \hspace*{-0.6in} - [B_{1k}, \Theta^+((C_1)_k)] \right) \right)^T  F.
\end{align}
\endth 

\section{Preservation of Commutation and Anti-Commutation Relations} \label{sec:section4}

In this section, the main results of the paper are provided. It is shown first under what conditions a state space model, described by the QSDE
\rref{eq:bilinear_system}, preserves the commutation and anticommutation relations of $SU(n)$. The employed procedure for this task is independent of the
physical realizability conditions given in Theorem \ref{th:physical_realizability}. Also, recall from Section \ref{sec:section3} that $G\in \C^{s\times s}$ can
be always regarded as $G \in \mathfrak{T}(\mathfrak{H})^{s\times s}$ since the identity operator in $\mathfrak{T}(\mathfrak{H})$ can be considered attached to
each component of $G$. The following Lemma plays a key role in the forthcoming theorems.

\begle \label{le:ZLemma} Let $G\in \mathfrak{T}(\mathfrak{H})^{s\times s}$ be antisymmetric and $x$ be a vector comprised by elements of a complete
linearly independent set of operators spanning $\mathfrak{T}(\mathfrak{H})$. If $G$ satisfies 
\begeq \label{eq:G_Theta_g}
G\Theta^-(x)+\Theta^-(x)G^T - \Theta^-(Gx) = 0,
\endeq
then there exist $g\in \mathfrak{T}(\mathfrak{H})$ such that
\begeq \label{eq:G_Theta_minus}
G=\Theta^-(g),
\endeq
where the unique $g$ is given by 
\begeq \label{eq:constructing_g_out_of_G}
g\triangleq  -\frac{1}{n} \colvec{3}{ \Tr(F_1 G) }{ \vdots }{ \Tr(F_s G) }. 
\endeq
Conversely, if \rref{eq:G_Theta_minus} holds, then \rref{eq:G_Theta_g} is true for any $x\in \mathfrak{T}(\mathfrak{H})^s$.
\endle 
\begpr Assume \rref{eq:G_Theta_g} holds. By applying $\vec$ to \rref{eq:G_Theta_g} and recalling that $F=(F_1,\cdots,F_s)^T$, it follows that
\begin{align*}
\lefteqn{\vec\left( G\Theta^-(x)+\Theta^-(x)G^T - \Theta^-(Gx) \right)} \\ 
& \hspace*{0.6in}= (I\otimes G)Fx + (G\otimes I)Fx - FGx. 
\end{align*}
Since by assumption the components of $x$ are linearly independent, it can be concluded that 
\begeq \label{eq:G_Theta_g_proof1}
(I\otimes G)F + (G\otimes I)F - FG = 0. 
\endeq
The application of identity $\rref{eq:kron_permutation3}$ in Lemma \ref{le:F_kronecker_properties} and \rref{eq:E^T_by_E} after multiplying on the left by $F^T$
gives
\begin{align*}
2\,\sum_{k=1}^s F_k G F_k + nG = 0.
\end{align*}
The equation for the $(i,j)$ component of $G$ is then
\begin{align} \label{eq:G_Theta_g_proof2}
2 \sum_{k,r,l=1}^s f_{kjr}f_{kil} G_{l r} - n G_{ij}=0. 
\end{align}
From \rref{eq:property_ff_tensor} and noting that $\Tr(AB)=\sum_{r,l=0}^s A_{rl}B_{lr}$ for any $A,B\in \re^{s\times s}$, \rref{eq:G_Theta_g_proof2} is
rewritten as
\begin{align} \label{eq:ZLemma1}
\nonumber \lefteqn{\hspace*{-0.21in}0 = \, 2 \sum_{k,r,l=1}^s f_{kjr}f_{kil} G_{l r} - n G_{ij}}  \\
\nonumber = & \, -2 \sum_{r,l=1}^s G_{l r} \sum_{k=1}^s (f_{jlk}f_{ikr} + f_{rlk}f_{ijk})  - n G_{ij}\\
\nonumber = & \, -2 \sum_{k,r,l=1}^s f_{jlk}f_{ikr} (G^T)_{r l} \\
\nonumber   & \, -2\sum_{k=1}^s f_{ijk} \sum_{r,l=0}^s f_{rlk} G_{l r}  - n G_{ij} \\
\nonumber = & \, -2 \sum_{k,r,l=0}^s f_{jlk}f_{ikr} (G^T)_{r l} \\
   & \, - 2 \sum_{k=1}^s f_{ijk} \Tr(F_k G)  -  n G_{ij}.
\end{align}
Also \rref{eq:G_Theta_g_proof2} implies that
\begin{align*}
2 \!\!\sum_{k,r,l=1}^s f_{jlk}f_{ikr} (G^T)_{r l} & = 2 \!\!\sum_{k,r,l=1}^s\!\! f_{kjl}(-f_{kir}) (-G_{rl }) =  n\, G_{ij}.  \\
\end{align*}
Substituting this on \rref{eq:ZLemma1} gives
\begdi
0= -n G_{ij} - 2 \sum_{k=1}^s f_{ijk} \Tr(F_k G) - n G_{ij},
\enddi
which means that
$G_{ij}=-\frac{1}{n} \sum_{k=1}^s f_{ijk} \Tr(F_k G)$. 
Then,
\begin{align*}
G_{ij} & = - \frac{1}{n} \sum_{k=1}^s f_{ijk} \Tr(F_k G)\\
& = \frac{1}{n} ( f_{ji1} , \cdots, f_{jis} ) \colvec{3}{ \Tr(F_1 G) }{ \vdots }{  \Tr(F_s G) } \\
& = - \frac{1}{n} ( (F_{j})_{i1}, \cdots, (F_{j})_{is} ) \colvec{3}{ \Tr(F_1 G) }{ \vdots }{  \Tr(F_s G) }.
\end{align*}
Therefore, 
\begin{align*}
G_j & = -\frac{1}{n} F_j^T  \colvec{3}{ \Tr(F_1 G) }{ \vdots }{  \Tr(F_s G) }.
\end{align*}
By defining $g$ as in \rref{eq:constructing_g_out_of_G}, one obtains that
\begdi
G = (F_1^T g, \cdots, F_s^T g ) = \Theta^-(g).
\enddi
Now, assuming that there exist a $g$ such that $G=\Theta^-(g)$, it then follows directly from \rref{eq:Theta__minus_theta_minus} that the left-hand side of
\rref{eq:G_Theta_g} satisfies
\begdi
\Theta^-(g)\Theta^-(x)-\Theta^-(x)\Theta^-(g)^T - \Theta^-(\Theta^-(g)x) = 0
\enddi
for an arbitrary $x\in \mathfrak{T}(\mathfrak{H})^{s}$, which concludes the proof.
\endpr

In order to be considered a quantum system, the system variables of \rref{eq:bilinear_system} must preserve the product rule
\rref{eq:product_Gell-mann_matrices} for all times. This amounts to preserve the commutation and anticommutation relations of $SU(n)$. Thus, the conditions
that
\rref{eq:bilinear_system} must satisfy are 
\begin{subequations}
    \begin{alignat}{2}
\label{eq:CCR_pauli}  d[x,x^T] & - 2{\pmb i} \Theta^-(dx) &\, = 0,\\
\label{eq:ACCR_pauli} d\{x,x^T\} & - 2 \Theta^+(dx) & = 0.
    \end{alignat}
\end{subequations}
Note by the linearity of $\Theta^-$ and $\Theta^+$ that
\begin{align*}
\Theta^-(dx)  = & \;\, \Theta^-(A_0) dt + \Theta^-(A x)  dt \\
& \,+ \Theta^-({B}_1 x) d\bar{W}_1+ \Theta^-({B}_2 x) d\bar{W}_2,\\
\Theta^+(dx)  = & \;\, \Theta^+(A_0) dt + \Theta^+(A x)  dt \\
& \,+ \Theta^+({B}_1 x) d\bar{W}_1+ \Theta^+({B}_2 x) d\bar{W}_2.
\end{align*}
A condition for system \rref{eq:bilinear_system} to satisfy \rref{eq:CCR_pauli} and \rref{eq:ACCR_pauli} is given in the next theorem.

\begth \label{th:preservation_commutation_relations} Let $x(0)\in \mathfrak{T}(\mathfrak{H})^s$ be comprised by the generators of $SU(n)$. Then, system
\rref{eq:bilinear_system} implies 
\begin{align*} 
[x(t),x(t)] & = 2{\pmb i} \Theta^-(x(t)), \\  
\{x(t),x(t)\} & = \frac{4}{n}I +2 \Theta^+(x(t))
\end{align*}
for all $t \ge 0$, if and only if
\begin{subequations}
    \begin{alignat}{2}
\label{subeqn:CCR_pauli_4} & {B}_i = \Theta^-(b_i) \\
\label{subeqn:CCR_pauli_1} & \sum_{k=1}^{n_w} {B}_{1k}{B}_{2k}^T - {B}_{2k}{B}_{1k}^T  = \frac{n}{2}\Theta^-(A_0) \\
\nonumber                  & A  = \Theta^-(a)  - \frac{1}{2}\sum_{k=1}^{n_w} \left(B_{1k}B_{1k}^T + B_{2k}B_{2k}^T\right) \\
\label{subeqn:CCR_pauli_5} &    \;\;\;\;\;\;\; + \frac{1}{2}\sum_{k=1}^{n_w} \left( B_{2k} \Theta^+(b_{1k}) - B_{1k} \Theta^+(b_{2k})\right).
    \end{alignat}
\end{subequations}
where $b_{ik}$ and $a$ are $s$-dimensional vectors as in \rref{eq:constructing_g_out_of_G} for $i=1,\hdots, s$ and $k=1,\hdots , n_w$.
\endth

\begpr
Without loss of generality, one can consider \rref{eq:bilinear_system} interacting with only one quadrature field. That is,
\begdi
dx = A_0\,dt + Ax\, dt + B_1 x\, d\bar{W}_1+B_2 x\, d\bar{W}_2.
\enddi
Using the It\^o table \rref{eq:Ito_table_quadrature}, it follows that
\begin{align*}
\lefteqn{\hspace*{-0.13in}d(xx^T) = (dx)x^T+x(dx)^T+(dx)(dx)^T}\\
	= \,& {} (A_0 x^T + x A_0^T)\,dt + (Axx^T+xx^TA^T)\,dt \\ 
	  & {}  + ({B}_1xx^T+xx^T{B}_1^T)\,d \bar{W}_1 +({B}_2xx^T+xx^T{B}_2^T)\,d\bar{W}_2  \\
	  & {}  + {B}_1x x^T {B}_1^T d\bar{W}_1 d\bar{W}_1+{B}_1x x^T {B}_2^T d\bar{W}_1 d\bar{W}_2 \\
	  & {}  + {B}_2x x^T {B}_1^T d\bar{W}_2 d\bar{W}_1 +{B}_2x x^T {B}_2^T d\bar{W}_2 d\bar{W}_2\\
	= \,& {} (A_0 x^T + x A_0^T)\,dt + (Axx^T+xx^TA^T)\,dt \\ 
	  & {}  + ({B}_1xx^T+xx^T{B}_1^T)\,d \bar{W}_1 +({B}_2xx^T+xx^T{B}_2^T)\,d\bar{W}_2  \\
	  & {}  + {B}_1x x^T {B}_1^T dt + \pmb{i}{B}_1x x^T {B}_2^T dt \\
	  & {}  - \pmb{i}{B}_2x x^T {B}_1^T dt + {B}_2x x^T {B}_2^T dt.
\end{align*}
Similarly, 
\begin{align*}
\lefteqn{\hspace*{-0.18in}\left(d(xx^T)\right)^T}\\
	= \,& {} (A_0 x^T + x A_0^T)\,dt + (A(xx^T)^T + (xx^T)^TA^T)\,dt \\ 
	  & {}  + ({B}_1(xx^T)^T + (xx^T)^T{B}_1^T)\,d \bar{W}_1 \\
	  & {}  +({B}_2(xx^T)^T + (xx^T)^T{B}_2^T)\,d\bar{W}_2  \\
	  & {}  + {B}_1 (xx^T)^T {B}_1^T dt + \pmb{i}{B}_1 (xx^T)^T {B}_2^T dt \\
	  & {}  - \pmb{i}{B}_2 (xx^T)^T {B}_1^T dt + {B}_2 (xx^T)^T {B}_2^T dt.
\end{align*}
Computing $d\left[x,x^T\right]$ then gives
\begin{align} \label{eq:preservation_CCR}
\nonumber d([x,x^T]) = & {\;} d(xx^T) - \left(d(xx^T)\right)^T \\
\nonumber            = & {\;} 2 \pmb{i} \left( \frac{2}{n}{B}_1{B}_2^T - \frac{2}{n}{B}_2{B}_1^T   \right)\,dt \\
\nonumber              & {\,} + 2\pmb{i}\!\left(\rule{0in}{0.2in}  A \Theta^-(x) +  \Theta^-(x) A^T \right)\, dt \\ 
\nonumber              & {\,} + 2\pmb{i}\!\left(\rule{0in}{0.2in} {B}_1 \Theta^-(x) {B}_1^T +  {B}_2 \Theta^-(x) {B}_2^T \right)\, dt \\ 
\nonumber              & {\,} + 2\pmb{i}\!\left( \rule{0in}{0.2in} {B}_1 \Theta^+(x) {B}_2^T -  {B}_2 \Theta^+(x) {B}_1^T \right)\, dt\\ 
\nonumber              & {\;} + 2\pmb{i} \left( {B}_1 \Theta^-(x) +  \Theta^-(x) {B}_1^T  \right) d\bar{W}_1 \\
                       & {\;} + 2\pmb{i} \left( {B}_2 \Theta^-(x) + \Theta^-(x) {B}_2^T  \right) d\bar{W}_2.
\end{align}
Replacing \rref{eq:preservation_CCR} into \rref{eq:CCR_pauli} amounts to 
\begin{align} \label{eq:preservation_CCR1}
\nonumber & {}\hspace{-0.1in} 2 \pmb{i}  \left( {B}_2{B}_1^T - {B}_1{B}_2^T - \Theta^-(A_0) +  A \Theta^-(x) +  \Theta^-(x) A^T \right.  \\
\nonumber & {}\hspace{-0.1in} \left. +{B}_1 \Theta^-(x) {B}_1^T +  {B}_2 \Theta^-(x) {B}_2^T - \Theta^-(Ax) \right) dt \\ 
\nonumber & {}\hspace{-0.1in}  + 2\pmb{i} \left( {B}_1 \Theta^-(x) +  \Theta^-(x) {B}_1^T -\Theta^-({B}_1x) \right) d\bar{W}_1  \\
	  & {}\hspace{-0.1in}  +  2\pmb{i} \left( {B}_2 \Theta^-(x) + \Theta^-(x) {B}_2^T -\Theta^-({B}_2x) \right) d\bar{W}_2=0.
\end{align}
From \cite[Proposition 27.3]{Parthasarathy_92}, one can also equate the integrands in \rref{eq:preservation_CCR1} to zero. Also, recall that that $x(0)$ is
represented by the complete orthonormal set constituted by the the generalized Gell-Mann matrices. This implies that any linear combination
$\sum_{k=0}^s a_i x_i(0) \not = 0$ unless $a_i=0$ for all $i$ and $a_i \in \C$. In addition, no linear combination of generalized Gell-Mann matrices generates
the identity. So, given that $x(0)\not = 0 $, any equation involving the system variables of the form $Ax=b$ ($A\in \C^{s\times s}$) implies $A$ and $b$ must be
identically $0$. Thus, the equations to be satisfied for preservation of commutation relations  are
\begin{subequations}
    \begin{alignat}{2}
\label{subeqn:CCR_pauli_3_x} & {B}_1 \Theta^-(x) + \Theta^-(x) {B}_1^T - \Theta^-({B}_1 x) = 0\\
\label{subeqn:CCR_pauli_4_x} & {B}_2 \Theta^-(x) + \Theta^-(x) {B}_2^T - \Theta^-({B}_2 x) = 0\\
\label{subeqn:CCR_pauli_5_x} & \frac{2}{n}{B}_1{B}_2^T -  \frac{2}{n}{B}_2{B}_1^T -  \Theta^-(A_0) = 0 \\
\nonumber                    & A\Theta^-(x)+ \Theta^-(x)A^T    +  {B}_1\Theta^-(x) {B}_1^T + {B}_2\Theta^-(x) {B}_2^T \\
\label{subeqn:CCR_pauli_2_x} & +  {B}_1\Theta^+(x) {B}_2^T - {B}_2\Theta^+(x) {B}_1^T-\Theta^-(Ax) =0.
    \end{alignat}
\end{subequations}
In the same vein for the anticommutator, one has that $d\left(\{x,x^T\}\right)$ is given by
\begin{align} \label{eq:preservation_CCR}
\nonumber \lefteqn{\hspace*{-0.08in}d(\{x,x^T\})
= d(xx^T) + \left(d(xx^T)\right)^T} \\
\nonumber              & \hspace*{-0.12in} = {\;} \frac{4}{n} \left( A + A^T + {B}_1{B}_1^T + {B}_2{B}_2^T   \right)\,dt \\
\nonumber              & \hspace*{-0.05in}{\,} + 2\left(\rule{0in}{0.2in}  A_0\,x^T  +  x \,A_0^T \right)\, dt \\ 
\nonumber              & \hspace*{-0.05in}{\,} + 2\left(\rule{0in}{0.2in}  A \Theta^+(x) +  \Theta^+(x) A^T \right)\, dt \\ 
\nonumber              & \hspace*{-0.05in}{\,} + 2\left(\rule{0in}{0.2in} {B}_1 \Theta^+(x) {B}_1^T +  {B}_2 \Theta^+(x) {B}_2^T \right)\, dt \\ 
\nonumber              & \hspace*{-0.05in}{\,} + 2\left( \rule{0in}{0.2in} {B}_1 \Theta^-(x) {B}_2^T -  {B}_2 \Theta^-(x) {B}_1^T \right)\, dt\\ 
\nonumber              & \hspace*{-0.05in}{\;} + 2 \left( \frac{2}{n}B_1 + \frac{2}{n} B_1^T + {B}_1 \Theta^+(x) +  \Theta^+(x) {B}_1^T  \right) d\bar{W}_1 \\
                       & \hspace*{-0.05in}{\;} + 2 \left( \frac{2}{n}B_2 + \frac{2}{n} B_2^T + {B}_2 \Theta^+(x) + \Theta^+(x) {B}_2^T  \right) d\bar{W}_2.
\end{align}
Using the same argument employed for the commutation relation of $x$, the integrands in \rref{eq:preservation_CCR1} are equal to zero. Thus, for preservation of
anticommutation relations one has that
\begin{subequations}
    \begin{alignat}{2}
\label{subeqn:ACCR_pauli_3_x} & {B}_1 \Theta^+(x) + \Theta^+(x) {B}_1^T - \Theta^+({B}_1 x) = 0\\
\label{subeqn:ACCR_pauli_4_x} & {B}_2 \Theta^+(x) + \Theta^+(x) {B}_2^T - \Theta^+({B}_2 x) = 0\\
\label{subeqn:ACCR_pauli_1_x} & A + A^T + {B}_1{B}_1^T + {B}_2{B}_2^T - \frac{n}{2}\Theta^+(A_0) = 0  \\
\label{subeqn:ACCR_pauli_0_x} & {B}_1 + {B}_1^T = {B}_2 + {B}_2^T  = 0  \\
\nonumber                     &  A_0\, x^T + x\,A_0^T  +  A\Theta^+(x)+ \Theta^+(x)A^T - \Theta^+(Ax)  \\
\nonumber                     & +  {B}_1\Theta^+(x) {B}_1^T + {B}_2\Theta^+(x) {B}_2^T \\
\label{subeqn:ACCR_pauli_2_x} & -  {B}_1\Theta^-(x) {B}_2^T + {B}_2\Theta^-(x) {B}_1^T = 0.
    \end{alignat}
\end{subequations}
Assume that \rref{subeqn:CCR_pauli_3_x}-\rref{subeqn:CCR_pauli_2_x} and \rref{subeqn:ACCR_pauli_3_x}-\rref{subeqn:ACCR_pauli_2_x} hold. Condition
\rref{subeqn:CCR_pauli_1} appear explicitly in \rref{subeqn:CCR_pauli_4_x}. From \rref{subeqn:CCR_pauli_3_x}, \rref{subeqn:CCR_pauli_4_x} and
\rref{subeqn:ACCR_pauli_0_x}, condition \rref{subeqn:CCR_pauli_4} is obtained by direct application of Lemma~\ref{le:ZLemma},  i.e., $B_i=\Theta^-(b_i)$ holds,
where $b_i$ is given in \rref{eq:constructing_g_out_of_G} for $i\in\{1,2\}$. In order to show \rref{subeqn:CCR_pauli_5}, $A_0$ need to be expressed in a more
convenient form. This is achieved by using \rref{eq:Theta__minus_theta_minus} and \rref{subeqn:CCR_pauli_4_x}. That is, 
\begin{align*}
\Theta^-(A_0) & = \frac{2}{n} \left(\Theta^-(b_2)\Theta^-(b_1) - \Theta^-(b_2)\Theta^-(b_1) \right) \\
              & = \frac{2}{n} \Theta^-\left(\Theta^-(b_2)b_1\right).
\end{align*}
Due to the linearity of $\Theta^-$, $A_0$ is then uniquely determined by $A_0=\frac{2}{n}\Theta^-(b_2)b_1$. Computing $\Theta^+(A_0)$ in terms of $b_1$ and
$b_2$ using \rref{eq:Theta__plus_theta_minus1} gives
\begin{align} \label{eq:A0_decomp}
\nonumber \Theta^+(A_0) = & \, \frac{2}{n}\Theta^+(\Theta^-(b_2)b_1)  \\
\nonumber               = & \, \frac{1}{n}\left(\Theta^+(b_2)\Theta^-(b_1) - \Theta^-(b_1)\Theta^+(b_2) \right. \\
                          & \;\;\;\;\, +\left. \Theta^-(b_2)\Theta^+(b_1) - \Theta^+(b_1)\Theta^-(b_2) \right).
\end{align}
Using \rref{subeqn:ACCR_pauli_1_x} and \rref{eq:A0_decomp}, one obtains the following expression
\begin{align} \label{eq:symmetric_Pbar}
\nonumber \lefteqn{\hspace*{-0.15in} \Theta^-(b_1)\Theta^-(b_1) + \Theta^-(b_2)\Theta^-(b_2)} \\
\nonumber        = &  \; A - \frac{1}{2} \left( \Theta^-(b_2)\Theta^+(b_1) - \Theta^-(b_1)\Theta^+(b_2) \right) \\
          &  \; +\left(A - \frac{1}{2} \left(   \Theta^-(b_2)\Theta^+(b_1) - \Theta^-(b_1)\Theta^+(b_2) \right)\right)^T\!\!\!.
\end{align}
Now let $\bar{P} \triangleq A - \frac{1}{2} \left(\Theta^-(b_2)\Theta^+(b_1) - \Theta^-(b_1)\Theta^+(b_2)\right)$. One can trivially decompose $\bar{P}$ in its
symmetric and antisymmetric parts. That is,
\begin{align*} 
\bar{P} = & \; \frac{1}{2}\left(\bar{P} + \bar{P}^T\right) + \frac{1}{2}\left(\bar{P}-\bar{P}^T\right). 
\end{align*}
Define $P=\frac{1}{2}\left(\bar{P}-\bar{P}^T\right)$. Note from \rref{eq:symmetric_Pbar} that 
\begdi 
\bar{P} + \bar{P}^T = \Theta^-(b_1)\Theta^-(b_1) + \Theta^-(b_2)\Theta^-(b_2). 
\enddi 
Hence, $P$ can be written as
\begin{align*}
P = & \; A - \frac{1}{2}\left( \Theta^-(b_1)\Theta^-(b_1) + \Theta^-(b_2)\Theta^-(b_2) \right) \\
    & \; - \frac{1}{2} \left(\Theta^-(b_2)\Theta^+(b_1) - \Theta^-(b_1)\Theta^+(b_2)\right)
\end{align*}
It is only left to show that $P$ can be written in terms of $\Theta^-$. The formula for $P$ allows to calculate $\Theta^-(Px)$ as 
\begin{align*}
\lefteqn{\hspace*{-0.1in}\Theta^-(P x)} \\
              = & \; \Theta^-(A x) - \frac{1}{2}\,\Theta^-\left( \Theta^-(b_1)\Theta^-(b_1)x + \Theta^-(b_2)\Theta^-(b_2)x \right.   \\
                & \;  + \left. \Theta^-(b_2)\Theta^+(b_1) x - \Theta^-(b_1)\Theta^+(b_2) x \right) \\
              = & \;  \Theta^-(A x) \\
                & \;  - \frac{1}{2}\, \left(\Theta^-\left( \Theta^-(b_1)\Theta^-(b_1)x\right) + \Theta^-\left(\Theta^-(b_2)\Theta^-(b_2)x \right) \right)  \\
                & \;  - \frac{1}{2}\, \left(\Theta^-\left(\Theta^-(b_2)\Theta^+(b_1) x\right) - \Theta^-\left(\Theta^-(b_1)\Theta^+(b_2) x \right) \right) \\
              = & \;  \Theta^-(A x) + \Theta^-(b_1)\Theta^-(x)\Theta^-(b_1) + \Theta^-(b_2)\Theta^-(x)\Theta^-(b_2) \\
                & \;  -\frac{1}{2}\left( \Theta^-(b_1)\Theta^-(b_1)\Theta^-(x) + \Theta^-(b_2)\Theta^-(b_2)\Theta^-(x) \right. \\
                & \;\;\;\;\,  + \Theta^-(x)\Theta^-(b_1)\Theta^-(b_1) + \Theta^-(x)\Theta^-(b_2)\Theta^-(b_2) \\
                & \;\;\;\;\,  + \Theta^-(b_2)\Theta^+(b_1)\Theta^-(x) + \Theta^-(b_2)\Theta^+(x)\Theta^-(b_1) \\
                & \;\;\;\;\,  - \Theta^-(b_1)\Theta^+(x)\Theta^-(b_2) - \Theta^-(x)\Theta^+(b_1)\Theta^-(b_2)\\
                & \;\;\;\;\,  - \Theta^-(b_1)\Theta^+(b_2)\Theta^-(x) - \Theta^-(b_1)\Theta^+(x)\Theta^-(b_2) \\
                & \;\;\;\;\,  \left.+\, \Theta^-(b_2)\Theta^+(x)\Theta^-(b_1) + \Theta^-(x)\Theta^+(b_2)\Theta^-(b_1) \right).                
\end{align*}
In the previous calculation, the antisymmetry of $\Theta^-$ was used as follows
\begin{align*}
\Theta^-\left(\Theta^+(b_2)x\right) & = \Theta^-(b_2)\Theta^+(x) + \Theta^-(x)\Theta^+(b_2) \\
                                    & = -\left(\Theta^-(b_2)\Theta^+(x) + \Theta^-(x)\Theta^+(b_2)\right)^T \\
                                    & =  \Theta^+(x)\Theta^-(b_2) + \Theta^+(b_2)\Theta^-(x).
\end{align*}
Also,
\begin{align*}
\lefteqn{\hspace*{-0.15in} P \Theta^-(x) + \Theta^-(x) P^T = A \Theta^-(x) + \Theta^-(x) A^T} \\
                & \;  - \frac{1}{2}\left( \Theta^-(b_1)\Theta^-(b_1)\Theta^-(x) + \Theta^-(b_2)\Theta^-(b_2)\Theta^-(x) \right. \\
                & \;\;\;\;\,  + \Theta^-(x)\Theta^-(b_1)\Theta^-(b_1) + \Theta^-(x)\Theta^-(b_2)\Theta^-(b_2) \\
                & \;\;\;\;\,  + \Theta^-(b_2)\Theta^+(b_1)\Theta^-(x) - \Theta^-(x)\Theta^+(b_1)\Theta^-(b_2)\\
                & \;\;\;\;\,  \left. - \, \Theta^-(b_1)\Theta^+(b_2)\Theta^-(x) + \Theta^-(x)\Theta^+(b_2)\Theta^-(b_1) \right) .                
\end{align*}
Hence,
\begin{align} \label{eq:equation_for_P}
\nonumber \lefteqn{ \hspace*{-0.15in} P \Theta^-(x) +  \Theta^-(x) P^T - \Theta^-(P x)} \\
\nonumber                = & \; A \Theta^-(x) +  \Theta^-(x) A^T - \Theta^-(A x) \\
\nonumber                  & \; -  \Theta^-(b_1)\Theta^-(x)\Theta^-(b_1) - \Theta^-(b_1)\Theta^-(x)\Theta^-(b_1) \\
                           & \; -  \Theta^-(b_1)\Theta^+(x)\Theta^-(b_2) - \Theta^-(b_2)\Theta^+(x)\Theta^-(b_1).
\end{align}
From \rref{subeqn:CCR_pauli_2_x}, it is thus evident that $P$ obeys
\begin{align*}
P \Theta^-(x) +  \Theta^-(x) P^T - \Theta^-(P x) = 0.
\end{align*}
Given that this equation is satisfied for all $t$, a direct application of Lemma \ref{le:ZLemma} when $t=0$ shows that $P = \Theta^-(a)$ for a vector $a$
as in \rref{eq:constructing_g_out_of_G}, which implies
\rref{subeqn:CCR_pauli_5}.

Conversely, assume that \rref{subeqn:CCR_pauli_4}-\rref{subeqn:CCR_pauli_5} are true. By Lemma \ref{le:ZLemma}, equations
\rref{subeqn:CCR_pauli_3_x}-\rref{subeqn:CCR_pauli_5_x} and \rref{subeqn:ACCR_pauli_1_x} are automatically satisfied. Equation \rref{subeqn:ACCR_pauli_0_x} it
is easily obtained by computing $A+A^T$ using \rref{subeqn:CCR_pauli_5}. Furthermore, equations \rref{subeqn:ACCR_pauli_3_x} and \rref{subeqn:ACCR_pauli_4_x}
are a direct result of \rref{eq:Theta__plus_theta_minus1}. The validity of \rref{subeqn:CCR_pauli_2_x} is shown by computing $A \Theta^-(x) + \Theta^-(x) A^T 
-\Theta^-(Ax)$. This calculation is analogous to the computation of \rref{eq:equation_for_P}, but in this case it is already known that $P=\Theta^-(a)$. Thus,
\rref{subeqn:CCR_pauli_2_x} is true given that \rref{subeqn:CCR_pauli_4}-\rref{subeqn:CCR_pauli_5} hold. To show that \rref{subeqn:ACCR_pauli_2_x} holds, one
first computes $\Theta^+(Ax)$ using \rref{subeqn:CCR_pauli_5} as follows
\begin{align*}
\lefteqn{\hspace{-0.05in}\Theta^+(Ax) = \Theta^+\left(\Theta^-(a)x\right) } \\
& \; + \frac{1}{2} \left( \Theta^-(b_1)\Theta^+\left(\Theta^-(b_1)x \right) - \Theta^+\left(\Theta^-(b_1)x \right) \Theta^-(b_1) \right. \\
& \;  \;\;\;\;\;\,\left. +\Theta^-(b_2)\Theta^+\left(\Theta^-(b_2)x \right) - \Theta^+\left(\Theta^-(b_2)x \right) \Theta^-(b_2) \right) \\
& \; + \frac{1}{2} \left( \Theta^-(b_2)\Theta^+\left(\Theta^+(b_1)x \right) - \Theta^+\left(\Theta^+(b_1)x \right) \Theta^-(b_2) \right. \\
& \;  \;\;\;\;\;\,\left. -\Theta^-(b_1)\Theta^+\left(\Theta^+(b_2)x \right) + \Theta^+\left(\Theta^+(b_2)x \right) \Theta^-(b_1) \right) .
\end{align*}
Now, from the symmetry of $\Theta^+$ and \rref{eq:thetaplus_thetaplus}, one observes that
\begin{align*}
\lefteqn{ \Theta^+\left(\Theta^+(b_i)x\right) } \\
                                    & = \Theta^+(b_2)\Theta^+(x) - \Theta^-(x)\Theta^-(b_2) - \frac{2}{n} b_i^T x I + \frac{2}{n} b_i x^T \\
                                    & = \left(\Theta^+(b_2)\Theta^+(x) - \Theta^-(x)\Theta^-(b_2) - \frac{2}{n} b_i^T x I + \frac{2}{n} b_i x^T\right)^T \\
                                    & = \Theta^+(x) \Theta^+(b_2) - \Theta^-(b_2) \Theta^-(x) - \frac{2}{n} b_i^T x I + \frac{2}{n} x b_i^T.
\end{align*}
Therefore,
\begin{align*}
\lefteqn{\hspace{-0.05in}\Theta^+(Ax) = \Theta^+\left(\Theta^-(a)x\right) } \\
& \; + \frac{1}{2} \left( \Theta^-(b_1)\Theta^+(x)\Theta^-(b_1) - \Theta^-(b_1)\Theta^-(b_1) \Theta^+(x) \right. \\
& \;  \;\;\;\;\;\, -\Theta^+(x)\Theta^-(b_1)\Theta^-(b_1) + \Theta^-(b_1)\Theta^+(x) \Theta^-(b_1)  \\
& \;  \;\;\;\;\;\, +\Theta^-(b_2)\Theta^+(x)\Theta^-(b_2) - \Theta^-(b_2)\Theta^-(b_2) \Theta^+(x)  \\
& \;  \;\;\;\;\;\, \left. -\,\Theta^+(x)\Theta^-(b_2)\Theta^-(b_2) + \Theta^-(b_2)\Theta^+(x) \Theta^-(b_2)  \right) \\
& \; + \frac{1}{2} \left(\rule{0in}{0.2in} \Theta^-(b_2)\Theta^+(b_1)\Theta^+(x) - \Theta^-(b_2)\Theta^-(x) \Theta^-(b_1) \right.  \\
& \;  \;\;\;\;\;\, - \frac{2}{n}  \Theta^-(b_2) b_1^T x + \frac{2}{n} \Theta^-(b_2) b_1 x^T \\
& \;  \;\;\;\;\;\, -\Theta^+(x)\Theta^+(b_1)\Theta^-(b_2) + \Theta^-(b_1) \Theta^-(x) \Theta^-(b_2)  \\
& \;  \;\;\;\;\;\, + \frac{2}{n}  b_1^T x\, \Theta^-(b_2) - \frac{2}{n}  x b_1^T \Theta^-(b_2)  \\
& \; \;\;\;\;\;\,  -\Theta^-(b_1)\Theta^+(b_2)\Theta^+(x) + \Theta^-(b_1)\Theta^-(x) \Theta^-(b_2) \\
& \;  \;\;\;\;\;\, + \frac{2}{n}\Theta^-(b_1)b_2^T x - \frac{2}{n} \Theta^-(b_1)b_2 x^T  \\
& \;  \;\;\;\;\;\, + \Theta^+(x) \Theta^+(b_2)\Theta^-(b_1) - \Theta^-(b_2) \Theta^-(x) \Theta^-(b_1)  \\
& \;  \;\;\;\;\;\, \left. - \, \frac{2}{n}  b_2^T x\,\Theta^-(b_1) + \frac{2}{n} x b_2^T \Theta^-(b_1) \right).
\end{align*}
Also,
\begin{align*}
\lefteqn{\hspace*{-0.15in} A \Theta^+(x) + \Theta^+(x) A^T } \\
              = & \;  \Theta^-(a) \Theta^+(x) + \Theta^+(x) \Theta^-(a)^T \\
                & \;  + \frac{1}{2} \left( \Theta^-(b_1)\Theta^-(b_1)\Theta^+(x) + \Theta^-(b_2)\Theta^-(b_2)\Theta^+(x) \right. \\
                & \;  \;\;\;\;\;\, + \Theta^+(x)\Theta^-(b_1)\Theta^-(b_1) + \Theta^+(x)\Theta^-(b_2)\Theta^-(b_2) \\
                & \;  \;\;\;\;\;\, + \Theta^-(b_2)\Theta^+(b_1)\Theta^+(x) - \Theta^+(x)\Theta^+(b_1)\Theta^-(b_2)\\
                & \;  \;\;\;\;\;\, \left.- \, \Theta^-(b_1)\Theta^+(b_2)\Theta^+(x) + \Theta^+(x)\Theta^+(b_2)\Theta^-(b_1) \right).                
\end{align*}
Since $\frac{2}{n}\Theta^-(b_2)b_1x^T = A_0 x^T$ and $\frac{2}{n}\Theta^-(b_2)b_1x^T = xA_0^T$, It then follows that
\begin{align} \label{eq:preser_ACCR_last}
\nonumber \lefteqn{\hspace*{-0.2in} A \Theta^-(x) + \Theta^-(x) A^T  -\Theta^+(Ax)} \\
\nonumber              = & \;  \Theta^-(a) \Theta^+(x) + \Theta^+(x) \Theta^-(a)^T - \Theta^+(\Theta^-(a)x)\\
\nonumber                & \;  + \Theta^-(b_1)\Theta^+(x)\Theta^-(b_1) + \Theta^-(b_2)\Theta^+(x)\Theta^-(b_2) \\
\nonumber                & \;  + \Theta^-(b_1)\Theta^-(x) \Theta^-(b_2) + \Theta^-(b_2)\Theta^-(x) \Theta^-(b_1) \\
                & \;  - \frac{1}{n} \left( A_0x^T + x A_0^T \right).
\end{align}
Identity \rref{eq:Theta__plus_theta_minus1} shows directly that 
\begdi
\Theta^-(a) \Theta^+(x) + \Theta^+(x) \Theta^-(a)^T - \Theta^+(\Theta^-(a)x)=0.
\enddi
Hence, \rref{eq:preser_ACCR_last} implies \rref{subeqn:ACCR_pauli_2_x}, which completes the proof.
\endpr

\begth
A physically realizable system satisfies the conditions of Theorem \ref{th:preservation_commutation_relations}. 
\endth
\begpr Consider, without loss of generality, that $n_w=1$. Define $C= \frac{1}{2}(C_1+\pmb{i}C_2)$. This means that $C$ is playing the role of $\Lambda$ in
Theorem \ref{th:physical_realizability_def}. By conditions $\rref{itm:theorem_physical_realizability2}$ and $\rref{itm:theorem_physical_realizability21}$ of
Theorem \ref{th:physical_realizability}, the fact that $B_1=\Theta^-(b_1)$ and $B_2=\Theta^-(b_2)$ is automatic by considering $b_1 = C_2$ and $b_2 = C_1$.
Next, one has from $\rref{itm:theorem_physical_realizability1}$ of Theorem \ref{th:physical_realizability} that 
\begdi 
A_0=\frac{1}{n}\left(\frac{1}{2}(B_1 - \pmb{i}B_2)(C_1+\pmb{i}C_2)^\dagger \right) = -\frac{2}{n}{\pmb i}\Theta^-(C)C^\dagger.
\enddi
Now, from \rref{eq:Theta__minus_theta_minus}, $\Theta^-(A_0)$ is
\begin{align*}
\frac{n}{2}\Theta^-(A_0) & = -\Theta^-(2{\pmb i}\Theta^-(C)C^\dagger) \\
	    & = -2i\left(\Theta^-(C)\Theta^-(C^\dagger)-\Theta^-(C^\dagger)\Theta^-(C)\right).
\end{align*}
Since $B_i$ is antisymmetric, it follows that 
\begdi
B_1B_2^T-B_2B_1^T = - 2\pmb{i} (\Theta^-(C)\Theta^-(C^\dagger)-\Theta^-(C^\dagger)\Theta^-(C)) .
\enddi
Therefore, \rref{subeqn:CCR_pauli_1} holds. Finally, since the considered system is physically realizable then adding \rref{subeqn:F} and its transpose gives 
\begin{align*}
\lefteqn{\hspace*{-0.5in} A + A^T + B_1 B_1^T + B_2 B_2^T} \\
                = & \; \Theta^+(b_2)\Theta^-(b_1) - \Theta^-(b_1)\Theta^+(b_2)  \\
                  & \; +  \Theta^-(b_2)\Theta^+(b_1) - \Theta^+(b_1)\Theta^-(b_2) \\
                = & \; \frac{n}{2} \, \Theta^+(A_0),
\end{align*} 
which concludes the proof.
\endpr

\noindent The implication of this theorem is that a state space model, as given by \rref{eq:bilinear_system} and \rref{eq:bilinear_system_output}, describes an
open $n$-level quantum system when it satisfies the physical realizability conditions since they also ensure preservation of commutation and anticommutation
relations of $SU(n)$. In addition, it is known that the physical realizability conditions can be employed to obtain the $(S,L,H)$ parametrization of that
quantum system.

\section{Conclusions} \label{sec:conclusions}

Conditions for preserving commutations and anticommutation relations have been provided for QSDE's of the form \rref{eq:bilinear_system}. These results used
explicitly the algebra of $SU(n)$, and in particular algebraic expressions involving the anomaly tensor $d$ had to be employed. Moreover, it was shown that the
physical realizability implies the preservation of those relations, and therefore under physical realizability conditions the system given by
\rref{eq:bilinear_system} and \rref{eq:bilinear_system_output} describes an open $n$-level quantum system.

\end{document}